\documentclass[10pt]{article}
\usepackage{mathrsfs}
\usepackage{amsthm}
\usepackage{amsmath}
\usepackage{graphicx}
\usepackage{color}
\usepackage{amsfonts}
\newtheorem{theorem}{Theorem}[section]

\newtheorem{lemma}[theorem]{Lemma}

\newtheorem{conjecture}[theorem]{Conjecture}

\numberwithin{equation}{section}
\normalsize

\begin{document}
\title{\textbf{Mean field limit for bias voter model on regular trees}}
\author{Xiaofeng Xue \thanks{\textbf{E-mail}: xuexiaofeng@ucas.ac.cn \textbf{Address}: School of Mathematical Sciences,
University of Chinese Academy of Sciences, Beijing 100049, China.}\\ University of Chinese Academy of Sciences}

\date{}
\maketitle

\noindent {\bf Abstract:}

In this paper we are concerned with bias voter models on trees and
lattices, where the vertex in state $0$ reconsiders its opinion at a
larger rate than that of the vertex in state $1$. For the process on
tree with product measure as initial distribution, we obtain a mean
field limit at each moment of the probability that a given vertex is
in state $1$ as the degree of the tree grows to infinity.
Furthermore, for our model on trees and lattices, we show that the
process converges weakly to the configuration where all the vertices
are in state $1$ when the rate at which a vertex in state $0$
reconsiders its opinion is sufficiently large. The approach of
graphical representation and the complete convergence theorem of
contact process are main tools for the proofs of our results.

\noindent {\bf Keywords:} bias voter model, mean field limit,
asymptotically independent, contact process, complete convergence
theorem.

\section{Introduction}\label{section one}
In this paper, we are concerned with the bias voter model on regular
trees. First we introduce the definition of this process on general
graphs. The bias voter model $\{\eta_t\}_{t\geq 0}$ on a graph $S$
is with state space $\{0,1\}^S$. That is to say, at each vertex
$x\in S$, there is a spin taking value from $\{0,1\}$. For any
configuration $\eta\in \{0,1\}^S$ and $x\in S$, we denote by
$\eta(x)$ the value of the spin at $x$.

At $t=0$, each spin takes a value from $\{0,1\}$ according to some
probability distribution. $\lambda>\theta>0$ are two constants. For
each vertex $x$ in state $0$ (resp. $1$), it waits for an
exponential time with rate $\lambda$ (resp. $\theta$) to choose a
neighbor $y$ uniformly. Then, the value of the spin at $x$ flips to
that of the spin at $y$. Therefore, $\{\eta_t\}_{t\geq 0}$ is a spin
system (see Chapter 3 of \cite{LIG1985}) with flip rates function
given by
\begin{equation}\label{equ 1.1 flip rate}
c(x,\eta)=
\begin{cases}
\frac{\lambda}{{\rm deg}(x)}\sum_{y:y\sim x}\eta(y)
\text{\quad if~}\eta(x)=0,\\
\frac{\theta}{{\rm deg}(x)}\sum_{y:y\sim x}[1-\eta(y)] \text{\quad
if~}\eta(x)=1
\end{cases}
\end{equation}
for any $(x,\eta)\in S\times \{0,1\}^S$, where we denote by $x\sim
y$ that $x$ and $y$ are neighbors and denote by ${\rm deg}(x)$ the
degree of $x$.

Intuitively, $0$ and $1$ are two opposite opinions of a topic.
Vertices in state $0$ (resp. $1$) are individuals holding the
opinion $0$ (resp. $1$). Each individual waits for an exponential
time to choose a neighbor randomly and take the neighbor's opinion
as its. The assumption that $\lambda>\theta$ can be considered as
that $0$ is a more controversial idea such that individuals holding
it prefer to reconsider their opinions.

In this paper, we assume that the initial distribution of the
process is the product measure with density $p$ for $p\in (0,1)$,
which is denoted by $\mu_p$. In other words,
\[
\mu_p(\eta:\eta(x)=1,\forall~x\in A)=p^{|A|}
\]
for any $A\subseteq S$.

We denote by $P_S^p$ the probability measure of the process
$\{\eta_t\}_{t\geq 0}$ on $S$ with initial distribution $\mu_p$. It
is natural to consider the estimation of $P^p_S(\eta_t(x)=1)$, which
is the probability that $x$ takes $1$ at the moment $t$. When $S$ is
the complete graph with $N$ vertices, which we denote by $C_N$, it
is easy to prove that $P^p_{C_N}(\eta_t(x)=1)$ satisfies the
following limit theorem such that
\begin{equation}\label{equ 1.2}
\lim_{N\rightarrow+\infty}P^p_{C_N}(\eta_t(x)=1)=\frac{pe^{(\lambda-\theta)t}}{1-p+pe^{(\lambda-\theta)t}}
\end{equation}
for any $t>0$.

Equation \eqref{equ 1.2} follows from a classic theory about density
dependent population processes constructed by Ethier and Kurtz (see
Chapter 11 of \cite{Ethier1986}). Let $N_t$ be the number of
vertices in state $1$ at the moment $t$, then
\begin{equation}\label{equ 1.3}
N_t\rightarrow
\begin{cases}
N_t+1 \text{~at rate~} \frac{\lambda}{N}(N-N_t)N_t, \\
N_t-1 \text{~at rate~} \frac{\theta}{N}(N-N_t)N_t.
\end{cases}
\end{equation}

When the initial distribution of $\{\eta_t\}_{t\geq 0}$ is $\mu_p$,
then according to Theorem 11.2.1 of \cite{Ethier1986} and \eqref{equ
1.3}, $N_t/N$ converges weakly to the solution $f(t,p)$ to the
following ODE
\[
\begin{cases}
\frac{d}{dt}f(t,p)=(\lambda-\theta)f(t,p)[1-f(t,p)],\\
f(0,p)=p
\end{cases}
\]
as $N$ grows to infinity. The mathematical expression of $f(t,p)$ is
exactly the right side of \eqref{equ 1.2}.

It is natural to ask whether \eqref{equ 1.2} holds for other
homogeneous graph $S$. We manage to prove that the answer is
positive when $S$ is a regular tree, which we denote by
$\mathbb{T}^N$. For mathematical details, see Section \ref{section
two}.

Since
\[\lim_{t\rightarrow+\infty}f(t,p)=1\]
when $\lambda>\theta$, it is natural to guess that $\eta_t$
converges weakly to $\delta_1$, the configuration where all vertices
are in state $1$, as $t$ grows to infinity.  For $S$ is a regular
tree or a lattice, we can prove that this guess is correct when
$\lambda/\theta$ is sufficiently large. For mathematical details,
see Section \ref{section two}.

When $\lambda=\theta$, our model degenerates to the classic voter
model introduced by Clifford and Sudbury in \cite{Cli1973}. In
\cite{Holley1975}, Holley and Liggett give an important dual
relationship between the classic voter model and the coalescent
random walks, which shows that any invariant measure of the classic
voter model is a convex combination of $\delta_1$ and $\delta_0$
when and only when two independent simple random walks on $S$ will
meet with probability one. More details can be found in Section 3.4
and Chapter 5 of \cite{LIG1985}. The classic voter model is also a
linear system (see Chapter 9 of \cite{LIG1985}), which makes the
process has some good properties, such as $\sum_{x\in S}\eta_t(x)$
is a martingale. When $\lambda>\theta$, the bias voter model can not
be described via a linear system and has no good duality properties,
which makes the classic approach to deal with voter models not be
valid.

\section{Main results}\label{section two}
In this section we give the main results of this paper. We denote by
$\mathbb{T}^N$ the regular tree with degree $N+1$. We obtain that
\eqref{equ 1.2} holds for the bias voter model on $\mathbb{T}^N$.

\begin{theorem}\label{theorem 2.1 main}
For any $t>0$ and $p\in (0,1)$,
\begin{equation}\label{equ 2.1}
\lim_{N\rightarrow+\infty}P^p_{\mathbb{T}^N}(\eta_t(x)=1)=\frac{pe^{(\lambda-\theta)t}}{1-p+pe^{(\lambda-\theta)t}}.
\end{equation}
\end{theorem}

Equation \eqref{equ 2.1} gives the limit of the probability that a
given vertex is in state $1$ at the moment $t$ as the degree of the
tree grows to infinity for any $t>0$. The limit function
\[
f(t,p)=\frac{pe^{(\lambda-\theta)t}}{1-p+pe^{(\lambda-\theta)t}}
\]
is usually called the mean-field limit. Please note that
$P_{\mathbb{T}^N}^p(\eta_t(x)=1)$ does not depend the choice of $x$
since $\mu_p$ is a translation invariant measure on
$\{0,1\}^{\mathbb{T}^N}$ and the flip rate function of
$\{\eta_t\}_{t\geq 0}$ given by \eqref{equ 1.1 flip rate} is also
translation invariant.

The proof of Theorem \ref{theorem 2.1 main} is in Section
\ref{section 3}. The core step of the proof is to show that
$\eta_t(x)$ and $\eta_t(y)$ are asymptotically independent for a
pair of neighbors $x$ and $y$ as the degree of the tree grows to
infinity.

It is natural to ask whether the counterpart of Theorem \ref{theorem
2.1 main} for the bias voter model on lattices $\mathbb{Z}^d,
d=1,2,\ldots$ holds. We guess the answer is positive but we have not
manage to prove that.

We denote by $\eta_t\Rightarrow \mu$ when the process
$\{\eta_t\}_{t\geq 0}$ converges weakly to a probability measure
$\mu$. That is to say, $\eta_t\Rightarrow \mu$ when and only when
\[
\lim_{t\rightarrow+\infty}Ef(\eta_t)=\int_{\{0,1\}^S}f(\eta)~\mu(d\eta)
\]
for any $f\in C(\{0,1\}^S)$. The mean-field limit $f(t,p)$ given by
\eqref{equ 2.1} satisfies that $f(t,p)\rightarrow 1$ as
$t\rightarrow+\infty$. So it is natural to guess that
$P_{\mathbb{T}^N}^p(\eta_t(x)=1)\rightarrow 1$ and therefore
$\eta_t\Rightarrow \delta_1$, the configuration where all the
vertices are in state $1$. The following two theorems show that this
guess holds for the bias voter model on trees and lattices when
$\lambda/\theta$ is sufficiently large.

\begin{theorem}\label{theorem 2.2}
For each $N\geq 2$, there is a constant $A(N)>0$ such that when
$\lambda/\theta>A(N)$, then
\[
\eta_t\Rightarrow\delta_1
\]
for the bias voter model $\{\eta_t\}_{t\geq 0}$ on $\mathbb{T}^N$
with initial distribution $\mu_p$ with $p\in (0,1)$. The sequence
$\{A(N)\}_{N\geq 2}$ satisfies that
\begin{equation}\label{equ 2.2}
\limsup_{N\rightarrow+\infty}\frac{A(N)}{\sqrt{N}}\leq 1.
\end{equation}
\end{theorem}

The main approach to prove Theorem \ref{theorem 2.2} is to compare
the bias voter model with a contact process on tree. The fact that
the strong survived contact process on tree satisfies the complete
convergence theorem is crucial for our proof. The limit theorem
\eqref{equ 2.2} of $A(N)$ follows from an important estimation of
the second critical value of the contact process on tree. For
mathematical details, see Section \ref{section 4}. For more about
the contact process on tree, see \cite{Pem1992} and Part 1 of
\cite{LIG1999}.

We denote by $\mathbb{Z}^d$ the lattice with degree $2d$. The
following theorem is a counterpart of Theorem \ref{theorem 2.2} for
the bias voter model on $\mathbb{Z}^d$.

\begin{theorem}\label{theorem 2.3}
For each $d\geq 1$ and the bias voter model $\{\eta_t\}_{t\geq 0}$
on $\mathbb{Z}^d$ with initial distribution $\mu_p$ with $p\in
(0,1)$, when $\lambda/\theta>4$, then
\[
\eta_t\Rightarrow\delta_1.
\]
\end{theorem}
The proof of Theorem \ref{theorem 2.3} is nearly the same analysis
as that of Theorem \ref{theorem 2.2}. The assumption
$\lambda/\theta>4$ relies on the fact that the critical value for
the contact process on $\mathbb{Z}^d$ is at most $2/d$. For more
details, see Section \ref{section 4}.

\section{Mean-field limit}\label{section 3}
In this section, we will prove Theorem \ref{theorem 2.1 main}. For
any $t>0$ and $p\in (0,1)$, we define
\[
f(t,p)=\frac{pe^{(\lambda-\theta)t}}{1-p+pe^{(\lambda-\theta)t}}.
\]
Since we are focused on the case where $S=\mathbb{T}^N$ in this
section, we rewrite $P^p_{\mathbb{T}^N}$ as $P^p_N$. First it is
easy to show that $f(t,p)$ is an upper bound of $P^p_N(\eta_t(x)=1)$
for any $t\geq 0$ and $N\geq 1$.

\begin{lemma}\label{lemma 3.1}
For any $t\geq 0$ and $N\geq 1$,
\begin{equation}\label{equ 3.1}
P_N^p(\eta_t(x)=1)\leq f(t,p).
\end{equation}
\end{lemma}

\proof

According to the flip rate function $c(x,\eta)$ of
$\{\eta_t\}_{t\geq 0}$ given by \eqref{equ 1.1 flip rate} and
Hille-Yosida Theorem,
\begin{align}\label{equ 3.2}
\frac{d}{dt}P^p_N(\eta_t(x)=1)=&\frac{\lambda}{N+1}\sum_{y:y\sim
x}P_N^p(\eta_t(x)=0,\eta_t(y)=1)\\
&-\frac{\theta}{N+1}\sum_{y:y\sim
x}P_N^p(\eta_t(x)=1,\eta_t(y)=0)\notag
\end{align}
for any $t>0$.

Since $\mu_p$ and $c(x,\eta)$ are translation invariant,
\[
P_N^p(\eta_t(x)=0,\eta_t(y)=1)=P_N^p(\eta_t(x)=1,\eta_t(y)=0)
\]
and does not rely on the choose of the neighbor $y$.

Therefore,
\begin{equation}\label{equ 3.3}
\frac{d}{dt}P^p_N(\eta_t(x)=1)=(\lambda-\theta)P_N^p(\eta_t(x)=1,\eta_t(y)=0),
\end{equation}
where $y$ is a fixed neighbor of $x$.

It is easy to check that the bias voter model is an attractive spin
system (see Section 3.2 of \cite{LIG1985}). Therefore, the two
events $\{\eta_t(x)=1\}$ and $\{\eta_t(y)=0\}$ are negative
correlated when the initial distribution is $\mu_p$ according to
Theorem 2.2.14 of \cite{LIG1985}.

As a result,
\begin{align*}
P_N^p(\eta_t(x)=1,\eta_t(y)=0)&\leq
P_N^p(\eta_t(x)=1)P_N^p(\eta_t(y)=0)\\
&=P_N^p(\eta_t(x)=1)[1-P_N^p(\eta_t(x)=1)]
\end{align*}
and hence
\begin{equation}\label{equ 3.4}
\frac{d}{dt}\Big[\log
\frac{P_N^p(\eta_t(x)=1)}{1-P_N^p(\eta_t(x)=1)}\Big]\leq
(\lambda-\theta)
\end{equation}
by \eqref{equ 3.3}.

According to \eqref{equ 3.4},
\begin{equation}\label{equ 3.5}
\log \frac{P_N^p(\eta_t(x)=1)}{1-P_N^p(\eta_t(x)=1)}-\log
\frac{p}{1-p}\leq (\lambda-\theta)t
\end{equation}
for any $t>0$.

Equation \eqref{equ 3.1} follows from \eqref{equ 3.5} directly.

\qed

To give a lower bound of $P_N^p(\eta_t(x)=1)$, we give another
description of the bias voter model $\{\eta_t\}_{t\geq 0}$. We are
inspired by the approach of graphical representation introduce by
Harris in \cite{Har1978} and the construction of stochastic
processes of spin systems with exchange dynamics introduced by
Durrett and Neuhauser in \cite{Durrett1994}. For any $x,y\in
\mathbb{T}^N, x\sim y$, we assume that $\{N_{(x,y)}(t):t\geq 0\}$ is
a Poisson process with rate $(\lambda+\theta)/(N+1)$. Please note
that we care the order of $x$ and $y$, so $N_{(x,y)}\neq N_{(y,x)}$.
We assume that all these Poisson processes are independent. At
$t=0$, each spin takes a value from $\{0,1\}$ according to the
distribution $\mu_p$. Then, the spin at $x$ may change its value
only at event times of $N_{(x,y)},y\sim x$. For any $t>0$, we define
\[
\eta_{t-}(x)=\lim_{s\uparrow t,s<t}\eta_s(x)
\]
as the value of the spin at $x$ at the moment just before $t$. For
any event time $s$ of $N_{(x,y)}$ for some $y\sim x$, we flip a coin
with head probability $\frac{\lambda}{\lambda+\theta}$ and tail
probability $\frac{\theta}{\lambda+\theta}$ at the moment $s$. If
$\eta_{s-}(x)=0$ (resp. $\eta_{s-}(x)=1$) and the result of the coin
flipping is head (resp. tail), then $\eta_s(x)=\eta_{s-}(y)$,
otherwise, $\eta_{s}(x)=\eta_{s-}(x)$. According to the properties
of exponential distribution, it is easy to check that the process
$\{\eta_t\}_{t\geq 0}$ evolving according to the rules above is a
bias voter model with flip rate function $c(x,\eta)$ given by
$\eqref{equ 1.1 flip rate}$.

For any $x\sim y$ and $T>0$, we define
\[
A_{x,y}(T)=\{N_{(x,y)}(T)=N_{(y,x)}(T)=0\}
\]
as the random event that the first event time of $N_{(x,y)}$ and
$N_{(y,x)}$ does not come before $T$.

Then,
\[
P_N^p(A_{x,y}(T))=e^{-\frac{2(\lambda+\theta)}{N+1}T}
\]
and hence
\begin{equation}\label{equ 3.6}
\lim_{N\rightarrow+\infty}P_N^p(A_{x,y}(T))=1
\end{equation}
for any $T>0$ and $p\in (0,1)$.

After all the prepared work, we can give the proof of Theorem
\ref{theorem 2.1 main} now.

\proof[Proof of Theorem \ref{theorem 2.1 main}]

According to \eqref{equ 3.3}, $P^p_N(\eta_t(x)=1)$ is increasing
with $t$, therefore by Lemma \ref{lemma 3.1},
\begin{equation}\label{equ 3.7}
p\leq P^p_N(\eta_t(x)=1)\leq f(T,p)<1
\end{equation}
for any $t\in [0,T]$ and $N\geq 1$.

For any $t\in [0,T]$, $x\sim y\in \mathbb{T}^N$.
\begin{align}\label{equ 3.8}
P_N^p(\eta_t(x)=1,\eta_t(y)=0)&\geq P_N^p\Big(\eta_t(x)=1,\eta_t(y)=0,A_{x,y}(T)\Big)\\
&=P_N^p\Big(\eta_t(x)=1,\eta_t(y)=0\Big|A_{x,y}(T)\Big)P_N^p\Big(A_{x,y}(T)\Big).\notag
\end{align}

For $x\sim y\in \mathbb{T}^N$, let
\[
C_y(x)=\{z\in \mathbb{T}^N: \text{~there is a path avoiding $y$ from
$x$ to $z$}\}.
\]
Conditioned on $A_{x,y}(T)$, $\{\eta_t(z):z\in C_y(x)\}_{t\leq T}$
and $\{\eta_t(w):w\in C_x(y)\}_{t\leq T}$ are independent when the
initial distribution is $\mu_p$, since vertices in $C_x(y)$ can not
exchange opinions with vertices in $C_y(x)$ before the moment $T$.

As a result,
\begin{align}\label{equ 3.9}
&P_N^p\Big(\eta_t(x)=1,
\eta_t(y)=0\Big|A_{x,y}(T)\Big)\notag\\
&=P_N^p\Big(\eta_t(x)=1\Big|A_{x,y}(T)\Big)P_N^p\Big(\eta_t(y)=0\Big|A_{x,y}(T)\Big)\\
&\geq
\frac{\Big[P_N^p\big(\eta_t(x)=1\big)-P_N^p\big(A^c_{x,y}(T)\big)\Big]
\Big[P_N^p\big(\eta_t(y)=0\big)-P_N^p\big(A^c_{x,y}(T)\big)\Big]}{\Big[P^p_N\big(A_{x,y}(T)\big)\Big]^2}\notag
\end{align}
for ant $t\in [0,T]$, where $A^c_{x,y}(T)$ is the complementary set
of $A_{x,y}(T)$.

By \eqref{equ 3.3}, \eqref{equ 3.8} and \eqref{equ 3.9}, for $t\in
[0,T]$,
\begin{equation}\label{equ 3.10}
\frac{d}{dt}P_N^p(\eta_t(x)=1)\geq(\lambda-\theta)P_N^p(\eta_t(x)=1)\big[1-P_N^p(\eta_t(x)=1)\big]G_t^p(x,y,N),
\end{equation}
where
\begin{align}\label{equ 3.11}
G_t^p(x,y,N)&=\frac{\Big[1-\frac{P_N^p\big(A^c_{x,y}(T)\big)}{P_N^p\big(\eta_t(x)=1\big)}\Big]\Big[1-\frac{P_N^p\big(A^c_{x,y}(T)\big)}
{P_N^p\big(\eta_t(x)=0\big)}\Big]}{P_N^p\big(A_{x,y}(T)\big)}\notag\\
&\geq
\frac{\Big[1-\frac{P_N^p\big(A^c_{x,y}(T)\big)}{p}\Big]\Big[1-\frac{P_N^p\big(A^c_{x,y}(T)\big)}
{1-f(T,p)}\Big]}{P_N^p\big(A_{x,y}(T)\big)}
\end{align}
according to \eqref{equ 3.7}.

By \eqref{equ 3.6} and \eqref{equ 3.11}, for any $\epsilon>0$ and
$T>0$, there exists $N(\epsilon,T)>0$ such that
\begin{equation}\label{equ 3.12}
G_t^p(x,y,N)\geq 1-\epsilon
\end{equation}
for any $N\geq N(\epsilon,T)$ and $t\in [0,T]$.

By \eqref{equ 3.10} and \eqref{equ 3.12},
\begin{equation}\label{equ 3.13}
\frac{d}{dt}P_N^p(\eta_t(x)=1)\geq
(\lambda-\theta)(1-\epsilon)P_N^p(\eta_t(x)=1)\big[1-P_N^p(\eta_t(x)=1)\big]
\end{equation}
for $N\geq N(\epsilon,T)$ and $t\in [0,T]$.

By \eqref{equ 3.13},
\[
\frac{d}{dt}\Big[\log
\frac{P_N^p(\eta_t(x)=1)}{1-P_N^p(\eta_t(x)=1)}\Big]\geq
(\lambda-\theta)(1-\epsilon)
\]
and hence
\begin{equation}\label{equ 3.14}
P_N^p(\eta_t(x)=1)\geq
\frac{pe^{(\lambda-\theta)(1-\epsilon)t}}{1-p+pe^{(\lambda-\theta)(1-\epsilon)t}}
\end{equation}
for $N\geq N(\epsilon, T)$ and $t\in [0,T]$.

By \eqref{equ 3.14},
\begin{equation}\label{equ 3.15}
\liminf_{N\rightarrow+\infty}P_N^p(\eta_t(x)=1)\geq
\frac{pe^{(\lambda-\theta)(1-\epsilon)t}}{1-p+pe^{(\lambda-\theta)(1-\epsilon)t}}
\end{equation}
for any $t>0$ and $\epsilon>0$.

Theorem \ref{theorem 2.1 main} follows from \eqref{equ 3.1} and
\eqref{equ 3.15} directly.

\qed

In the proof above, we show that $\eta_t(x)$ and $\eta_t(y)$ are
asymptotically independent as $N$ grows to infinity, since
$\eta_t(x)$ and $\eta_t(y)$ are independent conditioned on
$A_{x,y}(t)$ and the probability of $A_{x,y}(t)$ converges to $1$ as
$N$ grows to infinity. If we could show that $\eta_t(x)$ and
$\eta_t(y)$ are asymptotically independent for $x\sim y, x,y\in
\mathbb{Z}^d$ as $d$ grows to infinity, then we could extend Theorem
\ref{theorem 2.1 main} to the case where the bias voter model is on
the lattice. We will work on this problem as a further study.

\section{Weak convergence}\label{section 4}
In this section we will give the proofs of Theorem \ref{theorem 2.2}
and Theorem \ref{theorem 2.3}. After a scaling of the time, it is
easy to see that the limit behavior of $\{\eta_t\}_{t\geq 0}$ only
depends on $\lambda/\theta$, so in this section we assume that
$\theta=1$.

The proofs of Theorem \ref{theorem 2.2} and \ref{theorem 2.3} are
very similar, so we only give details of the proof of Theorem
\ref{theorem 2.2}. For Theorem \ref{theorem 2.3}, we only give a
sketch of the proof.

First we introduce the definition of the contact process
$\{\zeta_t\}_{t\geq 0}$ on $\mathbb{T}^N$. $\{\zeta_t\}_{t\geq 0}$
is a spin system with state space $\{0,1\}^{\mathbb{T}^N}$ and flip
rate function given by
\begin{equation}\label{equ 4.1 rate function of contact process}
c_1(x,\zeta)=
\begin{cases}
1 \text{\quad if~} \zeta(x)=1,\\
\frac{\lambda}{N+1}\sum_{y:y\sim x}\zeta(y) \text{\quad if~}
\zeta(x)=0
\end{cases}
\end{equation}
for any $(x,\zeta)\in \{0,1\}^{\mathbb{T}^N}$.

The contact process is first introduced by Harris in \cite{Har1974}.
Chapter 6 of \cite{LIG1985} and Part 1 of \cite{LIG1999} give a
detailed summary of main properties of the contact process.
Intuitively, the contact process describes the spread of an
infection disease. Vertices in state $1$ are infected while vertices
in state $0$ are healthy. An infected vertex waits for an
exponential time with rate $1$ to become healthy and a healthy
vertex is infected at a rate proportional to the number of infected
neighbors.

According to the basic coupling of spin systems (see Section 3.1 of
\cite{LIG1985}), we can also use $P_N^p$ to denote the probability
measure of the contact process $\{\zeta_t\}_{t\geq 0}$ on
$\mathbb{T}^N$ with initial distribution $\mu_p$. We write $P_N^p$
as $P_{N,\lambda}^p$ when we need to distinguish $\lambda$.

The following lemma shows that we can control the evolution of the
bias voter model $\{\eta_t\}_{t\geq 0}$ from below by the contact
process $\{\zeta_t\}_{t\geq 0}$, which is crucial for us to prove
Theorem \ref{theorem 2.2}.

\begin{lemma}\label{lemma 4.1}
Assume that $\{\eta_t\}_{t\geq 0}$ is the bias voter model with flip
rate function $c(x,\eta)$ given by \eqref{equ 1.1 flip rate} with
$\theta=1$ and $\{\zeta_t\}_{t\geq 0}$ is the contact process with
flip rate function $c_1(x,\zeta)$ given by \eqref{equ 4.1 rate
function of contact process}, then
\begin{equation}\label{equ 4.2}
P_{N,\lambda}^p(\eta_t(x)=0,\forall~x\in A)\leq
P_{N,\lambda}^p(\zeta_t(x)=0,\forall~x\in A)
\end{equation}
for any $A\subseteq \mathbb{T}^N$ and any $t\geq 0$.
\end{lemma}

\proof

For any $\eta,\zeta\in \{0,1\}^{\mathbb{T}^N}$, we write $\eta\geq
\zeta$ when and only when $\eta(x)\geq \zeta(x)$ for any $x\in
\mathbb{T}^N$.

By direct calculation, it is easy to check that
\begin{equation}\label{equ 4.3}
\begin{cases}
c(x,\eta)\geq c_1(x,\zeta) \text{\quad if~} \eta(x)=\zeta(x)=0,\\
c(x,\eta)\leq c_1(x,\zeta) \text{\quad if~} \eta(x)=\zeta(x)=1
\end{cases}
\end{equation}
for any $\eta\geq \zeta$.

By \eqref{equ 4.3} and Theorem 3.1.5 of \cite{LIG1985},
\begin{equation}\label{equ 4.4}
\eta_t\geq \zeta_t
\end{equation}
for any $t>0$ in the sense of coupling when $\eta_0$ and $\zeta_0$
have the same distribution $\mu_p$.

Equation \eqref{equ 4.2} follows from \eqref{equ 4.4} directly.

\qed

Now we introduce the second critical value of the contact process on
tree, $\lambda$ above which makes the complete convergence theorem
hold.

The contact process $\{\zeta_t\}_{t\geq 0}$ is an attractive spin
system (see Section 3.2 of \cite{LIG1985}), therefore
\[
P_{N,\lambda_1}^1(\exists~ t_n \uparrow +\infty,
\eta_{t_n}(x)=1,\forall~n\geq 1)\geq P_{N,\lambda_2}^1(\exists~ t_n
\uparrow +\infty, \eta_{t_n}(x)=1,\forall~n\geq 1)
\]
for $\lambda_1>\lambda_2$. As a result, it is reasonable to define
the following critical value for each $N\geq 2$,
\begin{equation}\label{equ 4.5}
A(N)=\sup\{\lambda: P_{N,\lambda}^1(\exists~ t_n \uparrow +\infty,
\eta_{t_n}(x)=1,\forall~n\geq 1)=0\}.
\end{equation}

$A(N)$ is called the second critical value of the contact process on
$\mathbb{T}^N$. When $\lambda>A(N)$, the contact process is called
strong survived. For more details, see Section 1.4 of
\cite{LIG1985}.

According to Theorem 1.4.65 of \cite{LIG1985},
\[
\limsup_{N\rightarrow+\infty}\sqrt{N}\frac{A(N)}{N+1}\leq 1,
\]
which is exactly equation \eqref{equ 2.2}.

The following lemma is a corollary of the complete convergence
theorem of strong survived contact process on tree. Please note that
we denote by $\delta_0$ the configuration in
$\{0,1\}^{\mathbb{T}^N}$ where all the vertices are in state $0$.

\begin{lemma}\label{lemma 4.2}
When $\lambda>A(N)$, then there is a probability measure
$\nu_\lambda$ on $\{0,1\}^{\mathbb{T}^N}$ such that
\begin{equation}\label{equ 4.6}
\nu_\lambda(\zeta:\zeta=\delta_0)=0
\end{equation}
and
\begin{equation}\label{equ 4.7}
\zeta_t\Rightarrow \nu_\lambda
\end{equation}
when $\zeta_0$ has probability distribution $\mu_p$ with $p\in
(0,1)$.
\end{lemma}

\proof

We denote by $\zeta_t^1$ the contact process with
$\zeta_0=\delta_1$. According to Theorem 3.2.3 and Theorem 6.1.6 of
\cite{LIG1985}, when $\lambda>A(N)$, there exists probability
measure $\nu_\lambda$ such that
\[
\zeta_t^1\Rightarrow \nu_\lambda
\]
and
\[
\nu_\lambda(\zeta:\zeta=\delta_0)=0.
\]

Let
\[
\tau=\inf\{t:\zeta_t=\delta_0\},
\]
then Theorem 1 of \cite{Zhang1996} shows that for any probability
measure $\mu$ on $\{0,1\}^{\mathbb{T}^N}$ and $\{\zeta_t\}_{t\geq
0}$ with initial distribution $\mu$,
\begin{equation}\label{equ 4.8}
\zeta_t\Rightarrow
P_\mu(\tau<+\infty)\delta_0+P_\mu(\tau=+\infty)\nu_\lambda
\end{equation}
when $\lambda>A(N)$.

When $\mu=\mu_p$ for $p\in (0,1)$, there are infinite many vertices
in state $1$ at $t=0$ with probbaility one and hence
\begin{equation}\label{equ 4.9}
P_N^p(\tau<+\infty)=0.
\end{equation}

Lemma \ref{lemma 4.2} follows from \eqref{equ 4.8} and \eqref{equ
4.9} directly.

\qed

Equation with form as \eqref{equ 4.8} is called the complete
convergence theorem, which shows that the process with any initial
distribution converges weakly to a convex combination of invariant
measures. In \cite{Bez1990}, Bezuidenhout and Grimmett show that the
complete convergence theorem holds for the contact process on
$\mathbb{Z}^d$. References \cite{Zhang1996} authored by Zhang and
\cite{Sal1998} authored by Salzano and Schonmann give two different
proofs of the complete convergence theorem of the strong survived
contact process on trees. In \cite{ChenX2009} and \cite{Yao2012},
Chen and Yao show that the complete convergence theorem holds for
contact process in a random environment on $Z^+\times Z^d$. In
\cite{Handjani1999}, Handjani shows that the complete convergence
theorem holds for the threshold-one voter model on $\mathbb{Z}^d$
such that the process with any initial distribution converges weakly
to a convex combination of three invariant measures.

By \eqref{equ 3.3}, $P_N^p(\eta_t(x)=1)$ is increasing with $t$, so
it is reasonable to define
\[
h(N,p)=\lim_{t\rightarrow+\infty}P_N^p(\eta_t(x)=1)
\]
for each $N\geq 1$ and $p\in (0,1)$. It is easy to see that
$\{\eta_t\}_{t\geq 0}$ with initial distribution $\mu_p$ converges
weakly to $\delta_1$ when and only when $h(N,p)=1$. The following
lemma shows that there is a subsequence of $\{\eta_t\}_{t\geq 0}$
converges weakly to a convex combination of $\delta_1$ and
$\delta_0$.

\begin{lemma}\label{lemma 4.3}
For $\{\eta_t\}_{t\geq 0}$ with initial distribution $\mu_p$ on
$\mathbb{T}^N$, there is a sequence $\{t_n\}_{n\geq 1}$ increasing
to infinity such that
\[
\eta_{t_n}\Rightarrow h(N,p)\delta_1+[1-h(N,p)]\delta_0.
\]
\end{lemma}

\proof

By \eqref{equ 3.3}, it is easy to see that
\[
\liminf_{t\rightarrow+\infty}P_N^p(\eta_t(x)=1,\eta_t(y)=0)=0
\]
for $x\sim y$. Otherwise, there would be $\alpha>0$ such that
$P_N^p(\eta_t(x)=1,\eta_t(y)=0)\geq \alpha$ for any $t>T_0$, where
$T_0$ is a sufficiently large number. Then, by \eqref{equ 3.3},
\[
P_N^p(\eta_t(x)=1)-P_N^p(\eta_{T_0}(x)=1)\geq
\alpha(t-T_0)\rightarrow +\infty
\]
as $t$ grows to infinity, which is contradictory.

Therefore, there exists sequence $\{t_n\}_{n\geq 1}$ increasing to
infinity such that
\begin{equation}\label{equ 4.10}
\lim_{n\rightarrow+\infty}P_N^p(\eta_{t_n}(x)=1,\eta_{t_n}(y)=0)=0
\end{equation}
for $x\sim y$.

Since $\{0,1\}^{\mathbb{T}^N}$ is a compact space, there is a
subsequence of $\{\eta_{t_n}\}_{n\geq 1}$ that converges weakly to a
probability measure on $\{0,1\}^{\mathbb{T}^N}$ according to the
Helly's selection theorem (see Theorem 3.2.6 of \cite{Dur2010}).
Without loss of generality, we can assume that
$\{\eta_{t_n}\}_{n\geq 1}$ is a convergent sequence itself.

We denote by $\varphi$ the limit distribution of $\eta_{t_n}$ as $n$
grows to infinity. Then, according to \eqref{equ 4.10},
\begin{equation}\label{equ 4.11}
\varphi(\eta(x)=1,\eta(y)=0)=0
\end{equation}
for any $x\sim y$.

By \eqref{equ 4.11},
\begin{align*}
\varphi(\eta:\eta\neq \delta_0,\delta_1)&=\varphi(\exists~x\sim y,\eta(x)\neq \eta(y))\\
&\leq \sum_{x\sim y}\varphi(\eta(x)=1,\eta(y)=0)+\sum_{x\sim
y}\varphi(\eta(x)=0,\eta(y)=1)=0.
\end{align*}

As a result, $\varphi$ is a convex combination of $\delta_1$ and
$\delta_0$. Since
\[
\varphi(\eta(x)=1)=\lim_{n\rightarrow+\infty}P_N^p(\eta_{t_n}(x)=1)=h(N,p)
\]
according to definition of $h$,
\[
\varphi=h(N,p)\delta_1+[1-h(N,p)]\delta_0.
\]

\qed

Finally we can give the proof of Theorem \ref{theorem 2.2}.

\proof[Proof of Theorem \ref{theorem 2.2}]

We only need to show that $h(N,p)=1$ when $\lambda>A(N)$. When
$\lambda>A(N)$, for any $\epsilon>0$, by \eqref{equ 4.6} in Lemma
\ref{lemma 4.2}, there exists finite subset $D$ of $\mathbb{T}^N$
such that
\begin{equation}\label{equ 4.12}
\nu_\lambda(\zeta(x)=0,\forall~x\in D)\leq \epsilon.
\end{equation}

By Lemma \ref{lemma 4.3}, there exists sequence $\{t_n\}_{n\geq 1}$
increasing to infinity such that
\begin{equation}\label{equ 4.13}
\lim_{n\rightarrow+\infty}P_{N,\lambda}^p(\eta_{t_n}(x)=0,\forall~x\in
D)=1-h(N,p).
\end{equation}

By Lemma \ref{lemma 4.2} and \eqref{equ 4.12},
\begin{equation}\label{equ 4.14}
\lim_{n\rightarrow+\infty}P_{N,\lambda}^p(\zeta_{t_n}(x)=0,\forall~x\in
D)=\nu_\lambda(\zeta(x)=0,\forall~x\in D)\leq \epsilon.
\end{equation}

By \eqref{equ 4.2}, \eqref{equ 4.13} and \eqref{equ 4.14},
\[
1-h(N,p)\leq \epsilon
\]
for any $\epsilon>0$.

As a result,
\[
h(N,p)=1
\]
when $\lambda>A(N)$ and the proof is complete.

\qed

At the end of this section, we give a sketch of the proof of Theorem
\ref{theorem 2.3}.

\proof[Proof of Theorem \ref{theorem 2.3}]

Let $\{\xi_t\}_{t\geq 0}$ be contact process on $\mathbb{Z}^d$ with
flip rate function given by
\[
c_2(x,\xi)=
\begin{cases}
1 \text{\quad if~}\xi(x)=1,\\
\frac{\lambda}{2d}\sum_{y:y\sim x}\xi(y)\text{\quad if~}\xi(x)=0
\end{cases}
\]
for any $(x,\xi)\in \mathbb{Z}^d\times \{0,1\}^{\mathbb{Z}^d}$

Let $\lambda(d)$ be the first critical value of $\{\xi_t\}_{t\geq
0}$, that is to say,
\[
\lambda(d)=\sup\{\lambda:\lim_{t\rightarrow+\infty}P_{\mathbb{Z}^d,\lambda}^1(\xi_t(x)=1)=0\}.
\]
It is shown in \cite{Bez1990} that the complete convergence theorem
holds for $\{\xi_t\}_{t\geq 0}$ when $\lambda>\lambda(d)$. Then,
according to a similar analysis with that in the proof of Theorem
\ref{theorem 2.2},
\[
\eta_t\Rightarrow \delta_1
\]
for the bias voter model $\{\eta_t\}_{t\geq 0}$ on $\mathbb{Z}^d$
with initial distribution $\mu_p$ with $p\in (0,1)$ when
$\lambda>\lambda(d)$.

According to Corollary 6.4.4 of \cite{LIG1985},
\[
\frac{\lambda(d)}{2d}\leq \frac{2}{d}
\]
and hence
\[
\lambda(d)\leq 4.
\]
Therefore, when $\lambda>4$, the bias voter model on $\mathbb{Z}^d$
with initial distribution $\mu_p$ with $p\in (0,1)$ converges weakly
to $\delta_1$.

\qed

\section{Two conjectures}\label{section 5}
In this section we propose two conjectures. The first one is about
the mean field limit of the bias voter model on lattices.

\begin{conjecture}\label{conjecture 5.1}
For $p\in (0,1)$,
\[
\lim_{d\rightarrow+\infty}P_{\mathbb{Z}^d}^p(\eta_t(x)=1)=\frac{pe^{(\lambda-\theta)t}}{1-p+pe^{(\lambda-\theta)t}}
\]
for any $t>0$.
\end{conjecture}
As we introduced in Section \ref{section 3}, the main difficulty to
prove Conjecture \ref{conjecture 5.1} is to show that $\eta_t(x)$
and $\eta_t(y)$ are asymptotically independent for $x\sim y, x,y\in
\mathbb{Z}^d$ as $d$ grows to infinity. Since there are infinite
many paths on the lattice from $x$ to $y$ avoiding the edge
connecting $x$ and $y$, our proof of Theorem \ref{theorem 2.1 main}
is not applicable for the case where the process is on the lattice.

The second conjecture is about the weak convergence of the process.
We guess that Theorem \ref{theorem 2.2} and Theorem \ref{theorem
2.3} hold under a generalized condition.

\begin{conjecture}\label{conjecture 5.2}
For any $\lambda>\theta$, $S=\mathbb{T}^d$ or $\mathbb{Z}^d$ with
$d\geq 1$,
\[
\eta_t\Rightarrow \delta_1
\]
for $\{\eta_t\}_{t\geq 0}$ on $S$ with initial distribution $\mu_p$
with $p\in (0,1)$.
\end{conjecture}

According to the proof of Theorem \ref{theorem 2.2}, the core step
to prove Conjecture \ref{conjecture 5.2} is to verify a claim that
the limit distribution of any convergent subsequence of
$\{\eta_t\}_{t\geq 0}$ puts no mass on $\delta_0$. However, for
$\lambda$ not large enough for the complete convergence theorem of
the contact process to hold, we have not find a way to prove this
claim yet.

We will work on this two conjectures as a further study and hope to
discuss with readers who are interested in them.

\quad

\textbf{Acknowledgments.} We are grateful to the financial support
from the National Natural Science Foundation of China with grant
number 11171342 and China Postdoctoral Science Foundation (No.
2015M571095).

{}
\end{document}